\documentclass[onefignum,onetabnum]{siamonline171218}
\usepackage[utf8]{inputenc}
\usepackage[T1]{fontenc}
\usepackage{microtype}
\usepackage{amsfonts}
\usepackage[UKenglish]{babel}
\usepackage{amssymb}
\usepackage{mathtools}
\usepackage{graphicx}
\usepackage[caption=false]{subfig}
\usepackage{makecell}
\usepackage{tabularx}
\usepackage{tikz}
\usepackage[acronym]{glossaries}
\usepackage[sort,compress]{cite}
\usepackage{amsopn}
\usepackage{enumitem}
\DeclareFontFamily{U}{matha}{\hyphenchar\font45}
\DeclareFontShape{U}{matha}{m}{n}{
      <5> <6> <7> <8> <9> <10> gen * matha
      <10.95> matha10 <12> <14.4> <17.28> <20.74> <24.88> matha12
      }{}
\DeclareSymbolFont{matha}{U}{matha}{m}{n}
\DeclareFontSubstitution{U}{matha}{m}{n}

\DeclareFontFamily{U}{mathx}{\hyphenchar\font45}
\DeclareFontShape{U}{mathx}{m}{n}{
      <5> <6> <7> <8> <9> <10>
      <10.95> <12> <14.4> <17.28> <20.74> <24.88>
      mathx10
      }{}
\DeclareSymbolFont{mathx}{U}{mathx}{m}{n}
\DeclareFontSubstitution{U}{mathx}{m}{n}
\DeclareMathDelimiter{\vvvert}{0}{matha}{"7E}{mathx}{"17}

\makeatletter
\providecommand*{\diff}{\@ifnextchar^{\DIfF}{\DIfF^{}}}

\def\DIfF^#1{%
  \mathop{\mathrm{\mathstrut d}}%
  \nolimits^{#1}\gobblespace}

\def\gobblespace{%
  \futurelet\diffarg\opspace}

\def\opspace{%
  \let\DiffSpace\!%
  \ifx\diffarg(%
  \let\DiffSpace\relax%
  \else%
  \ifx\diffarg\[%
  \let\DiffSpace\relax%
  \else%
  \ifx\diffarg\{%
  \let\DiffSpace\relax%
  \fi\fi\fi\DiffSpace}%
\makeatother

\DeclareMathOperator{\dive}{div}

\DeclarePairedDelimiterX{\abs}[1]{\lvert}{\rvert}{\ifblank{#1}{\:\cdot\:}{#1}}
\DeclarePairedDelimiterX{\norm}[1]{\lVert}{\rVert}{\ifblank{#1}{\:\cdot\:}{#1}}
\DeclarePairedDelimiterX{\scprod}[2]{\langle}{\rangle}{%
  \ifblank{#1}{\:\cdot\:}{#1}%
  \,,\mathopen{}%
  \ifblank{#2}{\:\cdot\:}{#2}%
}

\DeclarePairedDelimiterX{\set}[1]{\lbrace}{\rbrace}{%
  
  #1
}

\DeclarePairedDelimiterX{\setc}[2]{\lbrace}{\rbrace}{%
  #1\ifblank{#2}{}{\,\delimsize\vert\,\mathopen{}#2}}
\newcommand{\Natural}{\mathbb{N}}

\newcommand{\Real}{\mathbb{R}}

\ifpdf
  \DeclareGraphicsExtensions{.eps,.pdf,.png,.jpg}
\else
  \DeclareGraphicsExtensions{.eps}
\fi
\setlist[enumerate]{leftmargin=.5in}
\setlist[itemize]{leftmargin=.5in}

\newsiamremark{remark}{Remark}
\newsiamremark{hypothesis}{Hypothesis}
\crefname{hypothesis}{Hypothesis}{Hypotheses}
\newsiamthm{claim}{Claim}
\headers{%
  Theoretical Foundation of the Weighted Laplace Inpainting Problem}{%
  L. Hoeltgen and A. Kleefeld and I. Harris and M. Breuß}
\title{Theoretical Foundation of the Weighted Laplace Inpainting Problem}
\author{
  Laurent Hoeltgen\thanks{Institute for Mathematics, Brandenburg Technical University, Platz der
    Deutschen Einheit~1, 03046 Cottbus, Germany,
    (\email{hoeltgen@b-tu.de}, \email{breuss@b-tu.de}).}
  \and
  Andreas Kleefeld\thanks{Institute for Advanced Simulation, Forschungszentrum Jülich GmbH, Jülich
    Supercomputing Centre, Wilhelm-Johnen-Straße, 52425 Jülich, Germany
    (\email{a.kleefeld@fz-juelich.de}).}
  \and
  Isaac Harris\thanks{Department of Mathematics, Texas A\&{}M University, 621A Blocker Building,
    3368 TAMU College Station TX 77843, USA
    (\email{iharris@math.tamu.edu}).}
  \and
  Michael Breuß\footnotemark[1]
}
\ifpdf
\hypersetup{
  pdftitle={Theoretical Foundation of the Weighted Laplace Inpainting Problem},
  pdfauthor={L. Hoeltgen and A. Kleefeld and I. Harris and M. Breuß}
}
\fi
\begin{document}
\newacronym{PDE}{PDE}{partial differential equation}
\maketitle
\begin{abstract}
      Laplace interpolation is a popular approach in image inpainting using partial differential
      equations. The classic approach considers the Laplace equation with mixed boundary conditions.
      Recently a more general formulation has been proposed where the differential operator consists
      of a point-wise convex combination of the Laplacian and the known image data. We provide the
      first detailed analysis on existence and uniqueness of solutions for the arising mixed
      boundary value problem. Our approach considers the corresponding weak formulation and aims at
      using the Theorem of Lax-Milgram to assert the existence of a solution. To this end we have to
      resort to weighted Sobolev spaces. Our analysis shows that solutions do not exist
      unconditionally. The weights need some regularity and fulfil certain growth conditions. The
      results from this work complement findings which were previously only available for a discrete
      setup.
\end{abstract}
\section{Introduction}
\label{sec:introduction}
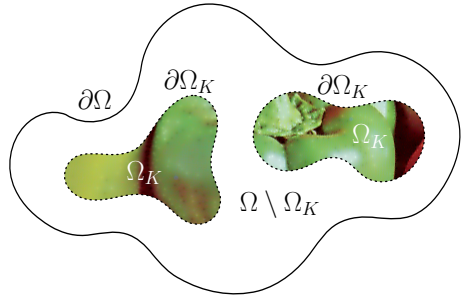
\begin{figure}
      \centering
      \begin{tikzpicture}[scale=0.6]
            \pgfdeclareimage[width=10cm]{background}{motivation2D-fill}
            \pgftext[base,left]{\pgfuseimage{background}};
            \pgftext[at=\pgfpoint{6cm}{2.0cm}]{\LARGE$\Omega\setminus{}\Omega_{K}$};
            \pgftext[at=\pgfpoint{8cm}{3.5cm}]{\textcolor{white}{\LARGE$\Omega_K$}};
            \pgftext[at=\pgfpoint{7.4cm}{4.6cm}]{\LARGE$\partial\Omega_K$};
            \pgftext[at=\pgfpoint{4.0cm}{4.7cm}]{\LARGE$\partial\Omega_K$};
            \pgftext[at=\pgfpoint{3cm}{2.7cm}]{\textcolor{white}{\LARGE$\Omega_K$}};
            \pgftext[at=\pgfpoint{1.9cm}{4.3cm}]{\LARGE$\partial\Omega$};
      \end{tikzpicture}
      \caption{Generic inpainting model as given in~\cref{eq:1} with known data image $f$ in
        $\Omega_K$. The task consists in recovering a reasonable reconstruction of the image $f$ in
        $\Omega\setminus{}\Omega_{K}$ by solving the \gls{PDE} in~\cref{eq:1}.}
      \label{fig:MBVP}
\end{figure}

Image inpainting deals with recovering lost image regions or structures by means of interpolation.
It is an ill-posed process; as soon as a part of the image is lost it cannot be recovered correctly
with absolute certainty unless the original image is completely known. The inpainting problem goes
back to the works of Masnou and Morel as well as Bertalmío and colleagues~\cite{MM1998a,BSCB2000},
although similar problems have been considered in other fields already before. There exist many
inpainting techniques, often based on interpolation algorithms, but \gls{PDE}[-based] approaches are
among the most successful ones, see e.g.~\cite{GM2014}. Among these, strategies based on the
Laplacian stand out~\cite{BW1989,SC2002a,MBWF2010,PHNH2016}. In that context, the elliptic mixed
boundary value problem
\begin{equation}
      \label{eq:1}
      \left\lbrace
            \begin{alignedat}{3}
                  -\Delta u &= 0
                  \quad &\text{in} &\quad &&\Omega\setminus\Omega_{K} \\
                  u &= f
                  &\text{in} &&&\partial\Omega_{K} \\
                  \partial_{n} u &= 0
                  &\text{in} &&&\partial\Omega\setminus{}\partial\Omega_{K}
            \end{alignedat}
      \right.
\end{equation}
is very popular. Here, $f$ represents known image data in a region $\Omega_{K}\subset\Omega$
(resp.~on the boundary $\partial\Omega_{K}$) of the whole image domain $\Omega$. Further,
$\partial_{n} u$ denotes the derivative in outer normal direction. An exemplary sketch of the layout
of the problem is given in~\cref{fig:MBVP}. Equations like~\cref{eq:1}, that involve different kinds
of boundary conditions, are commonly referred to as mixed boundary value problems and in rare cases
also as Zaremba's problem~\cite{Z1910}. Image inpainting based on~\cref{eq:1} appears under various
names in the literature: Laplace interpolation~\cite{PTVF2007}, harmonic
interpolation~\cite{S2015d}, or homogeneous diffusion inpainting~\cite{MBWF2010}. The latter name is
often used in combination with the steady state solution of the parabolic counterpart
of~\cref{eq:1}.\par
Applications of image inpainting are manifold and range from art restoration to image compression.
The earliest uses of \cref{eq:1} go back to Noma and Misulia (1959)~\cite{NM1959} and Crain
(1970)~\cite{C1970b} for generating topographic maps. Further applications include the works of
Bloor and Wilson (1989)~\cite{BW1989}, who studied partial differential equations for generating
blend surfaces. Finally, refer to~\cite{S2015d,HMHW2016} for a broad overview on \gls{PDE}[-based]
inpainting and the closely related problem of \gls{PDE}[-based] image compression.\par
\begin{figure}
      \centering
      \subfloat[Arbitrary data]{
        \label{fig:trui-demo-a}\includegraphics[width=0.24\linewidth]{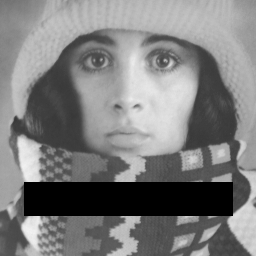}}
      \subfloat[Reconstruction]{
        \label{fig:trui-demo-b}\includegraphics[width=0.24\linewidth]{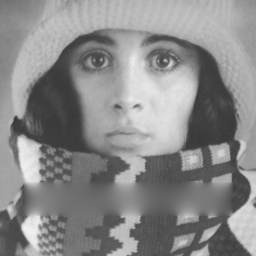}}\hfill
      \subfloat[Optimal data]{
        \label{fig:trui-demo-c}\includegraphics[width=0.24\linewidth]{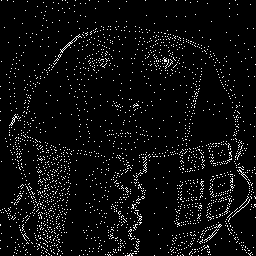}}
      \subfloat[Reconstruction]{
        \label{fig:trui-demo-d}\includegraphics[width=0.24\linewidth]{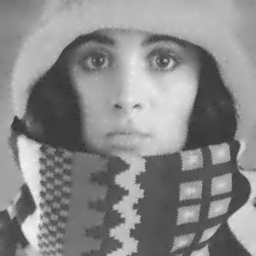}}
      \caption{(a) Image data with an arbitrary missing rectangular region (marked in black). (b)
        Corresponding reconstruction with \cref{eq:1}. The reconstruction suffers from blurring
        effects. (c) Remaining data (5\% of all pixels) with optimal reconstruction property.
        Missing data is black. (d) Corresponding reconstruction with \cref{eq:1}. The reconstruction
        is sharp although the Laplacian causes strong smoothing.}
      \label{fig:trui-demo}
\end{figure}

In the context of image reconstructions, \eqref{eq:1} is often favoured over other more complex
models due to its mathematically sound theory, even though the strong smoothing properties may yield
undesirable blurry reconstructions. Models based on anisotropic diffusion~\cite{GWWB2008,SWB2009} or
total variation~\cite{SC2002a} may be more powerful, but are much harder to grasp from a
mathematical point of view. In the context of image compression, the data $\Omega_{K}$ used for the
reconstruction can be freely chosen, since the complete image is known beforehand. The difficulty in
compressing an image with a \gls{PDE} lies in the fact that one has to optimise two contradicting
constraints. On the one hand, the size of the data $\Omega_{K}$ should be small to allow an
efficient coding, but on the other hand one wishes to have an accurate reconstruction from this
sparse amount of information, too. The optimal data depends on the choice of the differential
operator and the simplicity of the Laplacian offers many design choices for optimisation strategies
to find the best $\Omega_{K}$. Some of these approaches belong to the state-of-the-art methods in
\gls{PDE}[-based] image compression. We refer to~\cite{PHNH2016b} for a comparison of different
\gls{PDE}[-based] models and to \cite{GWWB2005,MHWT2012,HSW2013} for data optimisation strategies in
the compression context. \Cref{fig:trui-demo} demonstrates the potential of such a careful data
optimisation. \Cref{fig:trui-demo-a} and \cref{fig:trui-demo-b} show the reconstruction of an
arbitrary rectangle. The reconstruction is severely blurred and the texture of the scarf is almost
completely lost. On the other hand, \cref{fig:trui-demo-c} represents an optimised set of 5\% of the
data points with the corresponding colour values. \Cref{fig:trui-demo-d} depicts the corresponding
reconstruction. Although the reconstruction has a few artefacts, its overall quality is very
convincing.
As already mentioned, finding the best pixel-data is a very challenging task. Mainberger et
al.~\cite{MHWT2012} consider the combinatorial point of view of this task while Belhachmi and
colleagues~\cite{BBBW2009} approach the topic from the analytic side. Recently~\cite{HSW2013}, the
``hard'' boundary conditions in \cref{eq:1} have been replaced by softer weighting schemes that lead
to models such as
\begin{equation}
      \label{eq:2}
      \left\lbrace
            \begin{aligned}
                  c\left( u-f \right) + \left( 1-c \right) \left( -\Delta \right) u &= 0
                  & \text{in}\ & \Omega\setminus\Omega_{K} \\
                  \partial_{n} u &= 0
                  & \text{in}\ & \partial\Omega \\
            \end{aligned}
      \right.
\end{equation}
with a weighting function $c$. Optimising such a weighting function is notably simpler, at least in
discrete setups.\par
\Cref{eq:1} is well understood and there exist many results on existence, uniqueness, and regularity
of solutions, see~\cite{EG2004,CC2003} for a generic analysis and~\cite{SC2002a,CK2006} for a more
specific analysis in the inpainting context with Dirichlet boundary conditions only. Finite
difference discretisations of~\cref{eq:1} and~\cref{eq:2} have also been subject of several
investigations in the past. One can show that the discrete counterpart of~\cref{eq:1} admits a
unique solution as soon as the Dirichlet boundary set is non-empty~\cite{MBWF2010}. Similarly, the
discrete finite difference formulation of~\cref{eq:2} admits a unique solution if $c$ is positive in
at least one position~\cite{H2015b}.\par
An important question that arises in this context is what these discrete requirements relate to in
the continuous setting. If we consider for example the following model problem that one may extract
from the formulation of~\cref{eq:1},
\begin{equation}
      \label{eq:3}
      \left\lbrace
            \begin{alignedat}{3}
                  -\Delta u &= 0
                  \quad &\text{in} &\quad &&B_{1}\setminus B_{\varepsilon} \\
                  u &= 0
                  &\text{in} &&&\partial B_{1} \\
                  u &= 1
                  &\text{in} &&&\partial B_{\varepsilon}
            \end{alignedat}
      \right.
\end{equation}
where $B_{r}$ is a ball or radius $r$ with centre at the origin, then one can show that a smooth
solution exists for every $\varepsilon > 0$, but that no solution exists in the limiting case
$\varepsilon\to 0$. Indeed, the solution is given by
\begin{equation}
      \label{eq:4}
      u(x,y) = \frac{\ln\left( x^{2}+y^{2} \right)}{2\ln(\varepsilon)}
\end{equation}
Yet, the discrete formulation will admit a unique solution independently of the choice
of~$\varepsilon$. It suffices that the corresponding matrix is block irreducible. We refer
to~\cite{MBWF2010,H2015b} for a detailed discussion on the existence of solutions. To remedy the
situation for the continuous formulation in~\cref{eq:1}, the authors of~\cite{BBBW2009} have
required that the set $\Omega_{K}$ should have positive $\alpha$-capacity. The $\alpha$-capacity
($\alpha > 0$) of a subset $E\subset D$ of a smooth, bounded, and open set $D$ is given by
\begin{equation}
      \label{eq:43}
      \inf\;\setc*{\int_{D}\abs{\nabla u}^{2} + \alpha \abs*{u}^{2} \diff{x}}{u\in U_{E}}
\end{equation}
where $U_{E}$ is the set of all functions $u$ of the Sobolev space $H^{1}_{0}(D)$ such that $u
\geqslant{} 1$ almost everywhere in a neighbourhood of $E$. If $\Omega_{K}$ has positive
$\alpha$-capacity, then a solution of \cref{eq:1} exists in the Sobolev space
$H^{1}(\Omega)$~\cite{BBBW2009}. This requirement, that $\Omega_{K}$ must have positive capacity,
can be understood as requiring that image pixels are ``fat enough'' to allow a reconstruction. It
reconciles the continuous and discrete worlds and leads to a consistent theory on both sides. A
higher regularity than $H^{1}(\Omega)$ can be achieved for specific constellations of the boundary
data. A rather general theory is given in~\cite{F1949,M1970,AK1982}. The author of~\cite{M1955}
shows that a Hölder continuous solution exists if the data is regular enough. Finally,~\cite{B1994}
discusses the regularity of solutions on Lipschitz domains. Let us mention that the authors
of~\cite{CMS1998a} have also discussed this inability of the Laplacian to recover images from
isolated points and that they suggested absolutely minimising Lipschitz extensions as an
alternative.\par
The authors of this manuscript are not aware of any similar theory that would remedy the apparent
discrepancy between~\cref{eq:2} and its discrete counterpart. The discrete setup is almost always
solvable. On the other hand, solutions for the continuous model are only known for some special
cases such as $c$ being bounded between two positive constants in the interval $(0,1)$, or $c$ being
itself a constant~\cite{CC2003,EG2004}. For inpainting purposes it is important that $c$ may map to
the whole unit interval and even beyond. Regions with $c\equiv 1$ keep the data fixed and if $c$
exceeds the value 1, then contrast enhancing in the reconstruction is possible,
see~\cite{H2017e,HW2015a}.\par
Here, we attempt to bridge that gap between the discrete setup and the continuous model for the case
when $c$ maps to $[0,1]$. We show that a weak solution exists if certain assumptions on the weight
functions are met. Special interest will be paid to occurring requirements on $c$ and whether they
correspond to discrete counterparts. We aim at applying the theorem of Lax-Milgram in purpose-built
weighted Sobolev spaces. As a such, the contributed novelties of this manuscript are twofold. First
we complement the well-posedness study of \cref{eq:1} and $c > 1$, which has recently been discussed
in~\cite{HHBK2017} with the missing case where $c$ maps to $[0,1]$ and secondly, we introduce
weighted Sobolev spaces to the image processing community. These spaces bear a certain number of
interesting properties that can also be useful for other image analysis tasks.\par
In the next section we first derive the weak formulation corresponding to \cref{eq:1} and introduce
the weighted Sobolev spaces where the solution is sought. Then we will state the necessary
conditions on the weight function $c$ that must be fulfilled to assert the existence of a solution.
Finally we show that a unique solution exists.
%
\section{Inpainting with the weighted Laplacian}
\label{sec:inpa-with-lapl}
%
As already mentioned in the previous section, the classic formulation for \gls{PDE}[-based]
inpainting with the Laplacian reads
\begin{equation}
      \label{eq:5}
      \left\lbrace
            \begin{alignedat}{3}
                  -\Delta u &= 0
                  \quad &\text{in} &\quad &&\Omega\setminus\Omega_{K} \\
                  u &= f
                  &\text{in} &&&\partial\Omega_{K} \\
                  \partial_{n} u &= 0
                  &\text{in} &&&\partial\Omega\setminus{}\partial\Omega_{K}
            \end{alignedat}
      \right.
\end{equation}
Using the findings from~\cite{EG2004,HHBK2017}, it is easy to show that~\cref{eq:5} is well-posed
and that a unique weak solution exists in a subspace of $H^{1}(\Omega)$. With a weight function $c$
that maps from $\Omega$ to $\set{0,1}$, \cref{eq:5} can also be rewritten as
\begin{equation}
      \label{eq:6}
      \left\lbrace
            \begin{alignedat}{3}
                  c \left( u-f \right) + \left( 1-c \right)\left( -\Delta \right) u & = 0
                  \quad & \text{in} & \quad && \Omega \\
                  \partial_{n} u & = 0
                  & \text{in} &&& \partial\Omega \\
            \end{alignedat}
      \right.
\end{equation}
Interestingly, the latter formulation also makes sense if $c:\Omega\to\Real$, a fact which was first
exploited in~\cite{HSW2013}. If $c$ has binary values in the set $\set{0,1}$, then~\cref{eq:6} is
equivalent to~\cref{eq:5} with the Dirichlet boundary conditions specified by $f$ at those regions
where $c$ equals 1. \Cref{eq:6} can also be interpreted from a physical or chemical point of view.
We are in the presence of a stationary reaction-diffusion equation. The diffusive term
$(c-1)\Delta u$ is responsible for spreading the information generated by the reactive term
$c\, (u-f)$. The weight $c$ is responsible for the speed at which information is generated and
spread.\par
If $c$ is bounded between two non-negative numbers strictly smaller than one, then it follows
from~\cite{CC2003,HHBK2017} that a solution exists in $C^{2,\alpha}\left( \overline{\Omega}
\right)$. For inpainting purposes it is however important to allow $c(x)=1$ or even $c(x) > 1$. In
order to derive the weak formulation of \cref{eq:6} we follow the presentation in~\cite{HHBK2017},
where the setup in \cref{eq:6} with $c > 1$ was discussed by outlining its relationship to the
Helmholtz equation.\par
Let us now rewrite \eqref{eq:6} in a more suitable form. In a first step we explicitly set the
regions where $c\equiv 1$ apart.
\begin{equation}
      \label{eq:7}
      \left\lbrace
            \begin{aligned}
                  c\left( u-f \right) + \left( 1-c \right) \left( -\Delta \right) u &= 0
                  & \text{in}\ & \Omega\setminus\Omega_{K} \\
                  u &= f
                  & \text{in}\ & \partial\Omega_{K} \\
                  \partial_{n} u &= 0
                  & \text{in}\ & \partial\Omega\setminus{}\partial\Omega_{K} \\
            \end{aligned}
      \right.
\end{equation}
The previous reformulation implies that $c < 1$ almost everywhere in $\Omega\setminus\Omega_{K}$. A
small detail that will become important in the forthcoming discussions. If we further assume that
$c\in H^{1}(\Omega,[0,1])$, then we can apply the product rule and rewrite~\cref{eq:7} as
\begin{equation}
      \label{eq:8}
      \left\lbrace
            \begin{alignedat}{3}
                  -\dive \left( (1-c) \nabla u \right) - \nabla c \cdot \nabla u + c (u -f) &= 0
                  \quad & \text{in} & \quad & & \Omega\setminus\Omega_{K} \\
                  u &= f
                  & \text{in} &&& \partial\Omega_{K} \\
                  \partial_{n} u & = 0
                  & \text{in} &&& \partial\Omega\setminus{}\partial\Omega_{K} \\
            \end{alignedat}
      \right.
\end{equation}
In order to derive the weak formulation of~\cref{eq:8} within weighed Sobolev spaces let us first
remark that if $u$ solves \cref{eq:8}, then $v \coloneqq u - f$ solves
\begin{equation}
      \label{eq:15}
      \left\lbrace
            \begin{alignedat}{3}
                  -\dive \left( (1-c) \nabla v \right) - \nabla c \cdot \nabla v + c v &= g
                  \quad &\text{in} &\ &&\Omega\setminus\Omega_{K} \\
                  v &= 0
                  &\text{in} &&&\partial\Omega_{K} \\
                  \partial_{n} v &= h\
                  &\text{in} &&&\partial\Omega\setminus{}\partial\Omega_{K}
            \end{alignedat}
      \right.
\end{equation}
with $g \coloneqq (1-c) \Delta f$ and $h \coloneqq - \partial_{n} f$. For convenience of writing, we
will continue calling the sought solution of \cref{eq:15} $u$ and not $v$. Being able to solve
\cref{eq:15} is equivalent to being able to solve \cref{eq:8}. Yet, this change lets us reduce the
problem to the case with homogeneous Dirichlet boundary conditions. Deriving the associated weak
formulation is now straightforward. Multiplying with a test function $\varphi$ and integrating
\cref{eq:15} by parts implies that we must seek a function $u\in V$, which solves
\begin{equation}
      \label{eq:16}
      \underbrace{
        \int_{\mathrlap{\Omega\setminus{}\Omega_{K}}}\phantom{x}
        (1-c) \nabla u \cdot \nabla \varphi -
        (\nabla c \cdot \nabla u) \varphi +
        c u \varphi \diff{x}
      }_{\eqqcolon B^{c}(u,\varphi)}
      =
      \underbrace{\int_{\mathrlap{\Omega\setminus{}\Omega_{K}}}\phantom{x} g \varphi \diff{x} +
        \int_{\mathrlap{\partial\Omega\setminus{}\partial\Omega_{K}}}\phantom{x} h \varphi \diff{x}}_{\eqqcolon F(\varphi)}
      \quad\forall \varphi\in V
\end{equation}
Since $c$ maps to the unit interval, we are in the presence of a so called degenerate elliptic
equation~\cite{VG1969,SW2009} or sometimes also referred to as a \gls{PDE} with non-negative
characteristic form~\cite{OR1973}. Such \glspl{PDE} are characterised by the fact, that their
highest order term is allowed to vanish. This fact that the second order differential operator may
vanish locally requires a more sophisticated analysis. The key issue to approach this kind of
problems is to select the correct function space $V$ and to place certain necessary restrictions
onto $c$.\par
Let us briefly explain why additional restrictions on $c$, resp.\ the solution space $V$, are vital
to solve~\cref{eq:16}. The standard approach to show existence and uniqueness of a weak solution
consists in applying the Lax-Milgram Theorem~\cite{EG2004}. The crucial part will be the coercivity
of the bilinear form $B^{c}$ and the boundedness of $B^{c}$ and $F$. Obviously the boundedness of
$B^{c}$ and $F$ depends a lot on the choice of the space $V$ and $c$. To show coercivity of the
bilinear form, we must study the behaviour of
\begin{equation}
      \label{eq:10}
      \int_{\mathrlap{\Omega\setminus\Omega_{K}}}\phantom{x} (1-c) \abs*{\nabla u}^{2} -
      (\nabla c \cdot \nabla u) u + c u^{2} \diff{x}
\end{equation}
If $c$ is, for example, piecewise constant in~\cref{eq:10}, then the middle term vanishes almost
everywhere. Any function $u$ which is equal to 0 whenever $c$ is positive and equal to an arbitrary
constant when $c$ is 0 will force the bilinear form to be 0. Yet the norm of $u$ can be arbitrarily
large, which hinders us from showing coercivity. In order to prevent this situation, the following
assumptions seem reasonable:
\begin{enumerate}
\item The function $c$ should have a certain regularity, e.g.\ being continuous. Then arbitrary
      switching between regions where $c$ takes different constants is not possible anymore.
\item The space $V$ in which we seek solutions should fix its elements at certain boundaries. This
      would prevent solutions $u$ ``slipping'' away by adding constants that are invisible to the
      bilinear form.
\end{enumerate}
From these observations it becomes apparent that the function $c$ must go, at least partially, into
the definition of the space $V$. We consider such an approach in the following section by using
weighted Sobolev spaces and provide precise requirements that assert the well-posedness
of~\cref{eq:8}.
%
\subsection{Weighted Sobolev Spaces}
\label{sec:weight-sobol-spac}
%
Weighted Sobolev spaces have been studied intensively in the past. Their uses are manifold, but they
are most often found in the analysis of \glspl{PDE} with vanishing or singular diffusive term. The
works~\cite{VG1969,OR1973,SW2009,K1984,KS1987} give an almost complete overview of their usefulness.
For the sake of completeness, we shortly summarise how these spaces are set up.\par
In the following we denote by $W_{\Omega}$ the set of weight functions $\omega$, i.e.~$\omega$ is a
measurable and almost everywhere positive function in some domain $\Omega$. For $1 \leqslant{} p <
\infty$ and $\omega\in W_{\Omega}$ we define the corresponding weighted $L^{p}$ space as
\begin{equation}
      \label{eq:11}
      L^{p}(\Omega;\omega) \coloneqq \setc*{u\colon\Omega\to\Real}{
        \norm*{u}_{L^{p}(\Omega;\omega)} \coloneqq
        \left( \int_{\Omega} \abs*{u(x)}^{p}\omega(x)\diff{x} \right)^{\frac{1}{p}} < \infty}
\end{equation}
In a similar way as Sobolev spaces refine the Lebesgue spaces we can also refine our weighted
$L^{p}$ spaces by including the weak derivatives (defined in the usual sense) into the norm. In such
cases, different weights for different derivatives are also possible. For a given collection
$S_{k}\coloneqq\setc{\omega_{\alpha}\in W_{\Omega}}{\ \abs{\alpha} \leqslant k}$ of weight
functions, we denote by $W^{k,p}(\Omega; S_{k})$ the set of all functions $u$ defined on $\Omega$
whose (weak) derivatives $D^{\alpha}u$ of order $\abs*{\alpha} \leqslant{} k$ ($\alpha$ being a
multi-index) belong to $L^{p}\left( \Omega;\omega_{\alpha} \right)$. We can equip this vector space
with the norm
\begin{equation}
      \label{eq:12}
            \norm*{u}_{W^{k,p}(\Omega;S_{k})} \coloneqq
            \left(
                  \sum_{\abs*{\alpha}\leqslant{}k}
                  \int_{\Omega} \abs*{D^{\alpha}u(x)}^{p}\omega_{\alpha}(x)\diff{x}
            \right)^{\frac{1}{p}}
            =
            \left(
                  \sum_{\abs*{\alpha}\leqslant{}k}
                  \norm*{D^{\alpha}u}^{p}_{L^{p}(\Omega;\omega_{\alpha})}
            \right)^{\frac{1}{p}}
\end{equation}
One can show that the space $W^{k,p}(\Omega;S_{k})$ is a Banach space if $\omega_{\alpha}\in
L^{1}_{\text{loc}}(\Omega)$ and $\omega_{\alpha}^{\frac{-1}{p-1}}\in L^{1}_{\text{loc}}(\Omega)$ for
all $\abs{\alpha} \leqslant{} k$, see~\cite{KO1984,KO1986}. Note that this requires that all
derivatives up to the order $k$ must be attributed such a weight $\omega_{\alpha}$. However, one can
also show that $W^{k,p}(\Omega;\tilde{S}_{k})$ is still complete if $\tilde{S}_{k}\subsetneqq{}
S_{k}$ contains at least one weight $\omega_{\alpha}$ with $\abs{\alpha} = k$ and a weight for
$\abs{\alpha} = 0$, see~\cite{KO1983,KO1986a}.\par
We remark that for $p=2$ there is a canonical choice for a scalar product:
\begin{equation}
      \label{eq:14}
      \scprod*{u}{v}_{W^{k,2}(\Omega;S_{k})} \coloneqq
      \sum_{\abs*{\alpha}\leqslant{}k}
      \int_{\Omega} D^{\alpha}u(x)D^{\alpha}v(x)\omega_{\alpha}(x)\diff{x}
\end{equation}
Thus, with a suitable choice of weights we obtain a Hilbert space. If all the weight functions are
constant and equal to one, then our weighted spaces coincide with the usual definition of Sobolev
spaces. We refer to~\cite{K1984,KS1987} for a more complete listing of possible weighted Sobolev
space constructions. Finally, we remark that an alternative description of reasonable weight
functions can be given in terms of so called Muckenhoupt $A_{p}$ weights. We refer to~\cite{T2000b}
for more information.\par
By looking at~\cref{eq:8} it becomes apparent why these weighted Sobolev spaces are useful. The
function $c$ (resp.\ $1-c$) can be considered as a weight function and simply be integrated into the
space definition. This simplifies the proofs to show existence and uniqueness, since boundedness and
coercivity are easier to show and theorems such as Lax-Milgram can by applied in any Hilbert
space.\par
Our goal now will be to consider the corresponding weak formulation of~\cref{eq:16} in a suitable
weighted Sobolev space $V$. By applying the Theorem of Lax-Milgram in these spaces we will show the
existence and uniqueness of a weak solution of~\cref{eq:8}.\par
We make the following assumptions on our setup. These assumptions will hold throughout the whole
paper, unless mentioned otherwise. We remark that some of these assumptions can probably be
weakened, nevertheless they are not uncommon for image processing purposes and ease the discussion
on a few occasions.
\begin{enumerate}
\item $\Omega$ is an open, connected and bounded subset of $\Real^{2}$ with $\mathcal{C^{\infty}}$
      boundary $\partial\Omega$.\label{item:1}
\item $\Omega_{K}\subsetneqq\Omega$ is a closed subset of $\Omega$ with positive Lebesgue measure.
      It represents the known data locations used to recover the missing information on
      $\Omega\setminus{}\Omega_{K}$. The interpolation data is given by $f(\Omega_{K})$. The boundary
      $\partial\Omega_{K}$, is assumed to be $\mathcal{C}^{\infty}$, too. This set $\Omega_{K}$ is
      characterised by $c(x)\equiv 1$.\label{item:2}
\item $f\colon\Omega\to\Real$ is a $\mathcal{C}^{\infty}$ function representing the given
      image data to be interpolated by the underlying \gls{PDE}.\label{item:3}
\item The boundaries $\partial\Omega$ and $\partial\Omega_{K}$ do not intersect and neither of the
      boundaries $\partial\Omega$ or $\partial\Omega_{K}$ are empty.\label{item:4}
\item The function $c$ maps from $\Omega$ to the interval $[0,1]$, admits weak first order
      derivatives, and is an element of
      $L^{1}_{\text{loc}}\left( \Omega\setminus{}\Omega_{K} \right)$.\label{item:5}
\end{enumerate}
Let us briefly comment on these requirements. The first part of \cref{item:1} is trivially fulfilled
by images. Its second part is more restrictive. Assuming the boundary of $\Omega$ to be piecewise
$\mathcal{C^{\infty}}$ would be more realistic, but this would also reduce the regularity of the
solution. \Cref{item:2} and \cref{item:3} do not impose any severe restrictions for image processing
tasks. Images can always be rendered $\mathcal{C^{\infty}}$ by convolving them with a Gaussian.
\Cref{item:4} is necessary for technical reasons. If the Neumann and Dirichlet boundary conditions
meet each other, it is possible to generate setups that lead to contradicting requirements. Finally,
\cref{item:5} is necessary to assert the existence of our weighted Sobolev spaces.\par
The weights for our space definition should be chosen such that the bilinear form is equivalent to
the norm of our space. Often, the multiplicative factors of the individual derivatives in the
bilinear form offer themselves as viable choices for this task. In our case however, the function
$c$ may vanish locally. This prevents us from using $1-c$ and $c$ as weights to define a norm. They
only give us a seminorm structure. Such a situation is briefly described in~\cite{KO1983}. We mostly
follow that presentation and we propose the following correspondence between multi-indices
$\alpha\in \Natural_{0}^{2}$ and weights $\omega_{\alpha}$
\begin{equation}
      \label{eq:17}
      \omega_{\binom{0}{0}} \coloneqq 1,\quad
      \omega_{\binom{1}{0}} \coloneqq 1-c,\quad
      \omega_{\binom{0}{1}} \coloneqq 1-c
\end{equation}
This yields the scalar product and norm
\begin{subequations}
      \begin{align}
        \label{eq:18}
        \scprod{u}{v}_{V} &\coloneqq
                            \int_{\mathrlap{\Omega\setminus{}\Omega_{K}}}\phantom{x}
                            (1-c) \nabla u \cdot \nabla v + u v \diff{x} \\
        \label{eq:13}
        \norm*{u}_{V} &\coloneqq
                        \left( \int_{\mathrlap{\Omega\setminus{}\Omega_{K}}}\phantom{x}
                        (1-c) \abs*{\nabla u}^{2} + u^{2} \diff{x} \right)^{\frac{1}{2}}
      \end{align}
\end{subequations}
as well as the following definition for our space $V$:
\begin{equation}
      \label{eq:19}
      V \coloneqq
      \setc{
        \phi\in W^{1,2}\left(\Omega\setminus{}\Omega_{K};S_{c}\right)
      }{
        \phi\equiv 0\ \text{on}\ \partial\Omega_{K}
      }
\end{equation}
In addition, we define the following seminorm
\begin{equation}
      \label{eq:20}
      \left\vvvert u\right\vvvert_{V} \coloneqq
                  \left( \int_{\mathrlap{\Omega\setminus{}\Omega_{K}}}\phantom{x}
                  (1-c) \abs*{\nabla u}^{2} \diff{x} \right)^{\frac{1}{2}}
\end{equation}
Finally, following the presentation in~\cite{KS1987}, we note that the bilinear form $B^{c}$ in
\cref{eq:16} can be written compactly as a ternary quadratic form
\begin{equation}
      \label{eq:21}
      B^{c}(u,\varphi) = \sum_{\abs{\alpha},\abs{\beta} \leqslant{} 1}
      \int _{\mathrlap{\Omega\setminus{}\Omega_{K}}}\phantom{x}
      a_{\alpha,\beta} D^{\beta} u D^{\alpha} \varphi
      \diff{x}
\end{equation}
where $\alpha$, $\beta$ are multi-indices in $\Natural_{0}^{2}$. The weights $a_{\alpha,\beta}$ must
be set as follows to yield our model:
\begin{subequations}
      \label{eq:9}
      \begin{gather}
            \label{eq:22}
            a_{\binom{1}{0},\binom{1}{0}} = a_{\binom{0}{1},\binom{0}{1}} = 1 - c(x),\quad
            a_{\binom{0}{0},\binom{0}{0}} = c(x) \\
            \label{eq:33}
            a_{\binom{0}{0},\binom{1}{0}} = - \partial_{x} c(x),\quad
            a_{\binom{0}{0},\binom{0}{1}} = \partial_{y} c(x)
      \end{gather}
\end{subequations}
and $a_{\alpha,\beta} = 0$ for any other combination of multi-indices. In addition to the
previous assumptions, we assume further:
\begin{enumerate}
      \setcounter{enumi}{5}
\item There exists a constant $\kappa > 0$, such that for all $\abs{\alpha}$,
      $\abs{\beta} \leqslant{} 1$, $\alpha\neq\beta$,
      \begin{equation}
            \label{eq:23}
            \abs*{a_{\alpha,\beta}} \leqslant{} \kappa \sqrt{a_{\alpha,\alpha}a_{\beta,\beta}}
      \end{equation}
      almost everywhere in $\Omega\setminus\Omega_{K}$. For our choice in~\cref{eq:9}, this reduces to
      \begin{equation}
            \label{eq:24}
            \abs{\partial_{x} c} \leqslant{} \kappa \sqrt{c(1-c)},\quad
            \abs{\partial_{y} c} \leqslant{} \kappa \sqrt{c(1-c)}
      \end{equation}
      almost everywhere in $\Omega\setminus\Omega_{K}$.\label{item:6}
\item There exists a constant $\kappa' > 0$, such that for all real vectors $\xi\in\Real^{3}$ with
      entries $\xi_{\gamma}$ ($\gamma$ being a multi-index in $\mathbb{N}_{0}^{2}$ such that
      $\abs{\gamma} \leqslant 1$) we have
      \begin{equation}
            \label{eq:25}
            \sum_{\abs{\alpha},\abs{\beta}\leqslant{}1} a_{\alpha,\beta}\xi_{\alpha}\xi_{\beta}
            \geqslant{}
            \kappa'\sum_{\abs{\gamma}\leqslant{}1} a_{\gamma,\gamma}\xi_{\gamma}^{2}
      \end{equation}
      almost everywhere in $\Omega\setminus\Omega_{K}$. For our choice in~\cref{eq:9}, this reduces to
      \begin{subequations}
            \label{eq:44}
            \begin{align}
                  \label{eq:26}
                  &\begin{multlined}[][0.8\linewidth]
                        c\; \xi_{1}^{2} +
                        (1-c) \xi_{2}^{2} +
                        (1-c) \xi_{3}^{2}
                        - \partial_{x} c \; \xi_{1}\xi_{3} + \partial_{y} c\; \xi_{1}\xi_{2} \\
                        \geqslant{}
                        \kappa'
                        \left( (1-c) \xi_{3}^{2} + (1-c) \xi_{2}^{2} + c \xi_{1}^{2} \right)
                  \end{multlined} \\
                  \label{eq:27}
                  \Leftrightarrow{}\quad
                  &(\partial_{y} c) \xi_{1}\xi_{2} - (\partial_{x} c) \xi_{1}\xi_{3}
                  \geqslant{}
                  (\kappa'-1)
                  \left( (1-c) \xi_{3}^{2} + (1-c) \xi_{2}^{2} + c \xi_{1}^{2} \right)
            \end{align}
      \end{subequations}
      almost everywhere in $\Omega\setminus\Omega_{K}$.\label{item:7}
\end{enumerate}
\Cref{item:6} and \cref{item:7} are technical requirements that are necessary for the coercivity and
the boundedness of $B^{c}$. They cannot be avoided without substantial changes to the forthcoming
proofs. Let us remark, that~\cref{eq:25} can be deduced from~\cref{eq:23}, provided that
$\kappa < \frac{1}{2}$ holds. We refer to~\cite{KS1987} for a detailed proof. Equations~\cref{eq:24}
and~\cref{eq:27} enforce a certain well-behaviour on $c$, by restricting for example the growth
speed.\par
The following findings are a direct consequence of the foregoing results.
\begin{proposition}
      \label{thm:1}
      The bilinear form $B^{c}$ from~\cref{eq:21} is continuous.
      \begin{proof}
            By using~\cref{eq:23} and the Hölder inequality we obtain.
            \begin{equation}
                  \label{eq:28}
                  \begin{split}
                        \abs*{B^{c}(u,\varphi)} &\leqslant
                        \sum_{\abs{\alpha},\abs{\beta} \leqslant 1} \int _{\mathrlap{\Omega\setminus{}\Omega_{K}}}\phantom{x}
                        \abs{a_{\alpha,\beta}} \abs{D^{\beta} u} \abs{D^{\alpha} \varphi} \diff{x} \\
                        &\leqslant{} \max\set{\kappa, 1} \sum_{\abs{\alpha},\abs{\beta} \leqslant 1}
                        \int _{\mathrlap{\Omega\setminus{}\Omega_{K}}}\phantom{x}
                        \abs{D^{\beta} u} \sqrt{\abs{a_{\beta,\beta}}}
                        \abs{D^{\alpha} \varphi} \sqrt{\abs{a_{\alpha,\alpha}}} \diff{x} \\
                        &\leqslant{} K \norm{D^{\beta} u}_{V} \norm{D^{\alpha} \varphi}_{V}\\
                  \end{split}
            \end{equation}
            where $K$ is some positive constant. We emphasise that the last estimate requires
            $c \leqslant{} 1$ almost everywhere to be valid.
      \end{proof}
\end{proposition}
\begin{proposition}
      There exists a constant $\kappa' >0$ such that the bilinear form $B^{c}$ from~\cref{eq:21}
      satisfies the estimate $B^{c}(u,u) \geqslant{} \kappa'\left\vvvert u\right\vvvert_{V}^{2}$.
      \begin{proof}
            We replace $\xi_{\alpha}$ by $D^{\alpha} u$ and $\xi_{\beta}$ by $D^{\beta} u$ in
            \cref{eq:25}. Integrating the resulting inequality over $\Omega\setminus{}\Omega_{K}$
            yields
            \begin{equation}
                  \label{eq:29}
                  B^{c}(u,u)
                  =
                  \sum_{\abs{\alpha},\abs{\beta}\leqslant{}1}
                  \int_{\mathrlap{\Omega\setminus{}\Omega_{K}}}\phantom{x}
                  a_{\alpha,\beta}D^{\alpha} u D^{\beta} u
                  \diff{x}
                  \geqslant{}
                  \kappa'\sum_{\abs{\gamma}\leqslant{}1} \int_{\mathrlap{\Omega\setminus{}\Omega_{K}}}\phantom{x}
                  a_{\gamma,\gamma} \left( D^{\gamma} u \right)^{2}\diff{x}
                  \geqslant{}
                  \kappa'\left\vvvert u\right\vvvert_{V}^{2}
            \end{equation}
      \end{proof}
\end{proposition}
To complete the proof of the coercivity of the bilinear form $B^{c}$ we need a Friedrichs-like
estimate of the form $\norm{u}_V \leqslant{} K \left\vvvert u\right\vvvert_{V}$ with some positive
constant $K$. The particular formulation and preliminaries that we need can be found
in~\cite{WSZ2006} as Theorem~2.3. We repeat it here verbatim for the sake of completeness but refer
to its source for a detailed proof.\par
In the following theorem we denote by $W_{c}(X)$ the subset of weights on the space $X$ which are
bounded from above and below by positive constants on each compact subset $Q\subset X$, i.e.~we only
allow our weights to degenerate at the boundary of the domain. The next theorem also considers a
constant $A$ which is defined as follows. For an arbitrary domain $X$ we assume that we can
write
\begin{equation}
      \label{eq:30}
      X = \bigcup_{k=1}^{\infty} X_{k}
\end{equation}
where $(X_{k})_{k}$ is a sequence of bounded domains whose boundary can be locally described by
functions satisfying a Lipschitz condition and where $X_{k}\subset\overline{X}_{k}\subset X_{k+1}$
holds for each $k$. Finally, let $X^{k} \coloneqq X\setminus{}X_{k}$ and define
\begin{equation}
      \label{eq:31}
      A_{k} = \sup_{\norm{u}_{W^{k,p}(X; S_{k})} \leqslant{} 1} \norm{u}_{L^{p}(X^{k}; w_{0})}
\end{equation}
where $w_{0} \in S_{k}$ is the weight that corresponds to $\abs*{\alpha}=0$. We define additionally
$A \coloneqq \lim_{k\to\infty} A_{k}$. Obviously $A\in [0,1]$ always holds. This number $A$ is also
the ball measure of non-compactness of the embedding $W^{k,p}(X; S_{k}) \to L^{p}(X; w_{0})$,
see~\cite{WSZ2006,EO1993}. One can interpret the number $A$ as the distance from the embedding
operator to the next closest compact operator from $W^{k,p}(X; S_{k})$ into $L^{p}(X; w_{0})$. Also,
the numbers $A_{k}$ can be understood as indicators on how much ``weight'' is put onto the function
along the boundary. $A_{k} < 1$ means that there is at least some weight on the derivatives or
inside the domain. Note that in our setup~\cref{eq:31} simplifies to
\begin{equation}
      \label{eq:32}
      A_{k} = \sup_{
        \norm{u}_{W^{1,2}(\Omega\setminus{}\Omega_{K};S_{k})} \leqslant{} 1} \norm{u}_{L^{2}(X^{k})}
\end{equation}
where $X^{k}$ is the complement of a set $ X_{k}\subset\Omega\setminus{}\Omega_{K}$ and where
$S_{k}$ is the set of weights from~\cref{eq:17}.\par
For the following theorem it is important that $A < 1$, i.e.~the weight is not completely
concentrated on the boundary. Let us remark that this requirement is in accordance with the discrete
theory established in~\cite{MBWF2010,H2015b}. In the discrete setting, there should be at least one
position with positive weight in the interior of the domain.\par
Let us emphasise that for our task at hand, such a construction with the requirement that $A < 1$ is
an additional regularity assumption on our image data $f$ and the mask function $c$. Indeed, part of
the boundary of the domain that we consider is fixed where $c\equiv 1$. Since the $\Omega_{k}$ need
boundaries that can be described locally by functions that fulfil a Lipschitz condition, this
requirement carries over to the function $c$.\par
As already mentioned, the next theorem is a almost verbatim copy of Theorem 2.3 in~\cite{WSZ2006}.
\begin{theorem}
      \label{thm:2}
      Suppose $1 \leqslant{} p < \infty$ and $S_{k} \subset W_{c}(X)$. Let $\ell$ be a functional on
      $W^{k,p}\left(X; S_{k}\right)$ with the following properties.
      \begin{enumerate}
      \item $\ell$ is continuous on $W^{k,p}\left( X; S_{k}\right)$
      \item $\ell(\lambda u) = \lambda \ell(u)$ for all $\lambda > 0$ and all
            $u\in W^{k,p}\left(X, S_{k}\right)$.
      \item If $u\in P_{k-1}\cap W^{k,p}\left(X; S_{k}\right)$ ($P_{k-1}$ being the set of all
            polynomials on $\Real^{n}$ of degree less than $k$) and $\ell(u) = 0$, then $u = 0$.
      \end{enumerate}
      Let $A < 1$. Then there is a constant $\kappa_{0}$ such that
      \begin{equation}
            \label{eq:34}
            \int_{X}\abs*{u}^{p} w_{0} \diff{x} \leqslant{}
            \kappa_{0} \left(
                  \abs*{\ell(u)}^{p} + \sum_{\abs{\alpha}=k} \norm*{D^{\alpha} u}_{L^{p}(X;w_{\alpha})}^{p}
            \right)
      \end{equation}
      Here, $w_{0}$ is the weight that corresponds to $\abs{\alpha}=0$.
\end{theorem}
The previous theorem can be seen as a generalisation to weighted spaces of a well-known theorem for
constructing equivalent norms out of seminorms in regular Sobolev spaces. See Theorem~7.3.12
in~\cite{AH2009}. \Cref{eq:34} can also be considered as a higher dimensional generalisation of the
Hardy inequality. We refer to~\cite{OK1990} for an extensive treatise on this inequality.\par
We now use~\cref{thm:2} with $p=2$, $k=1$, $n=2$, $w_{0} \equiv 1$, $w_{\alpha} = 1 - c$ for all
$\alpha$ and
\begin{equation}
      \label{eq:35}
      \ell(u) = \int_{\mathrlap{\partial\Omega_{K}}}\phantom{x} u \diff{x}
\end{equation}
With these choices we obtain the Friedrichs' inequality in our space $V$:
\begin{equation}
      \label{eq:36}
      \norm{u}^{2}_{L^{2}(\Omega\setminus{}\Omega_{K})} \leqslant{}
      \kappa_{0}
      \left\vvvert u\right\vvvert_{V}^{2}
\end{equation}
\Cref{eq:36} is the final key building block in showing the existence and uniqueness of a solution
of our \gls{PDE}. It allows us to show the coercivity of our bilinear form.
\begin{proposition}
      \label{thm:3}
      If~\cref{eq:36} holds, i.e.~the requirements of~\cref{thm:2} are fulfilled for the choice of
      $\ell$ from~\cref{eq:35} and for our selection of weights for our space $V$, then the bilinear
      form $B^{c}$ from~\cref{eq:21} is coercive.
      \begin{proof}
            Equation~\cref{eq:36} immediately implies that
            $\norm{u}_{V}^{2} \leqslant{} (1+\kappa_{0})\left\vvvert u\right\vvvert_{V}^{2}$. In
            combination with~\cref{eq:29} it follows that
            \begin{equation}
                  \label{eq:37}
                  B^{c}(u,u) \geqslant{}
                  \kappa'\left\vvvert u\right\vvvert_{V}^{2}
                  \geqslant{}
                  \frac{\kappa'}{1+\kappa_{0}} \norm{u}_{V}^{2}
            \end{equation}
      \end{proof}
\end{proposition}
\Cref{thm:3} completes the analysis of our bilinear form $B^{c}$. It remains to show that the
right-hand side of our weak formulation is continuous if we want to apply the Theorem of Lax-Milgram.
This final step is done in the following proposition.
\begin{proposition}
      \label{thm:4}
      The linear operator $F$ from~\cref{eq:16} is continuous, provided that $g$, $\Delta f$, and
$\frac{\nabla f}{\sqrt{1-c}}$ are in $L^{2}(\Omega\setminus{}\Omega_{K})$.
      \begin{proof}
            We remark that $\varphi\in V$ is 0 along $\partial\Omega_{K}$, and thus we can extend
            the boundary integral over that part. Using the Hölder inequality and Green's first
            identity, we obtain
            \begin{equation}
                  \label{eq:38}
                  \begin{split}
                  \abs{F(\varphi)} &\leqslant{}
                  \int_{\mathrlap{\Omega\setminus{}\Omega_{K}}}\phantom{x} \abs{g} \abs{\varphi} \diff{x} +
                  \abs*{\int_{\mathrlap{\partial\Omega\setminus{}\partial\Omega_{K}}}\phantom{x} h\varphi \diff{x}} \\
                  &\leqslant{}
                  \norm{g}_{L^{2}(\Omega\setminus{}\Omega_{K})}\norm{\varphi}_{L^{2}(\Omega\setminus{}\Omega_{K})} +
                  \abs*{\int_{\mathrlap{\Omega\setminus{}\Omega_{K}}}\phantom{x} \Delta f \varphi + \nabla f \cdot \nabla\varphi \diff{x}}\\
                  &\leqslant{}
                  \norm{g}_{L^{2}(\Omega\setminus{}\Omega_{K})}\norm{\varphi}_{V} +
                  \norm{\Delta f}_{L^{2}(\Omega\setminus{}\Omega_{K})}\norm{\varphi}_{V} +
                  \abs*{\int_{\mathrlap{\Omega\setminus{}\Omega_{K}}}\phantom{x} \nabla f \cdot \nabla\varphi \diff{x}}\\
            \end{split}
      \end{equation}
      The last integral can be estimated as follows
      \begin{equation}
            \label{eq:39}
            \begin{split}
                  \abs*{\int_{\mathrlap{\Omega\setminus{}\Omega_{K}}}\phantom{x} \nabla f \cdot \nabla\varphi \diff{x}}
                  &=
                  \abs*{\int_{\Omega\setminus{}\Omega_{K}} \frac{\nabla f}{\sqrt{1-c}}\sqrt{1-c}\,\nabla\varphi \diff{x}} \\
                  &\leqslant{}
                  \norm*{\frac{\nabla f}{\sqrt{1-c}}}_{L^{2}(\Omega\setminus{}\Omega_{K})}\norm*{\nabla\varphi}_{L^{2}(\Omega\setminus{}\Omega_{K};1-c)} \\
                  &\leqslant{}
                  \norm*{\frac{\nabla f}{\sqrt{1-c}}}_{L^{2}(\Omega\setminus{}\Omega_{K})}\norm*{\varphi}_{V}
            \end{split}
      \end{equation}
      Therefore, it follows that
      \begin{equation}
            \label{eq:40}
            \abs{F(\varphi)} \leqslant{}
            \left( \norm{g}_{L^{2}(\Omega\setminus{}\Omega_{K})} +
            \norm{\Delta f}_{L^{2}(\Omega\setminus{}\Omega_{K})} +
            \norm*{\frac{\nabla f}{\sqrt{1-c}}}_{L^{2}(\Omega\setminus{}\Omega_{K})} \right)
            \norm{\varphi}_{V}
      \end{equation}
      Thus, $F$ is a bounded linear functional.
      \end{proof}
\end{proposition}
We can now combine our results to prove our main result.
\begin{theorem}
      The weak formulation~\cref{eq:16} of the mixed boundary value problem~\cref{eq:15} has a
      unique solution in the space $V$. In addition, we know that
      \begin{equation}
            \label{eq:41}
            \norm*{u}_{V} \leqslant{} \frac{1 + \kappa_{0}}{\kappa'} \norm*{F}_{V^{*}}
      \end{equation}
      where $ \frac{\kappa'}{1 + \kappa_{0}}$ is the constant from~\cref{eq:37}. Here, $V^{*}$
      denotes the dual space of $V$.
      \begin{proof}
            From \cref{thm:1} and \cref{thm:3} it follows that our bilinear form $B^{c}$ is bounded
            and coercive. \Cref{thm:4} shows that the corresponding right-hand side $F$ is bounded,
            too. Therefore, from the Theorem of Lax-Milgram (see~\cite{EG2004}) it follows that
            there exists a unique $u\in V$ such that $B^{c}(u,\varphi) = F(\varphi)$ holds for all
            $\varphi\in V$. In addition, this $u$ fulfils
            $\norm{u} \leqslant{} \frac{1 + \kappa_{0}}{\kappa'} \norm*{F}_{V^{*}}$
      \end{proof}
\end{theorem}
Our weighted Sobolev space $V$ might be unsuited for other applications as it is very problem
specific. Having an embedding from $V$ to some standard Sobolev space $W^{k,2}$ would be very useful
in view of the many embedding theorems for these latter spaces which can be used to show a higher
regularity of the solution. It would also help in comparing solutions obtained when $c$ maps to the
unit interval but does not reach 0 or 1. In that case, the solutions live in $W^{1,2}$. By
construction of our space $V$, we immediately obtain $u\in L^{2}(\Omega\setminus\Omega_{K})$.
However, this result does not even acknowledge the existence of the weak derivatives. Any claims
beyond that are difficult to do. There exist a certain number of results concerning the embedding of
weighted Sobolev spaces into other spaces, however, their assumptions are often very abstract or
quite restrictive, e.g. all weights must be identical. We refer the reader to the
works~\cite{GO1988,GO1989,GO1991} for a detailed analysis on this topic.
%
\subsection{What happens if \texorpdfstring{$c\geqslant{}1$}{c greater equal 1}?}
\label{sec:what-if-c}
%
Let us shortly discuss the consequences of $c$ exceeding its upper limit 1. Similar conclusions can
also be drawn for the case $c \leqslant{} 0$, however, this latter situation usually does not occur
in practice.\par
There are no restrictions on $c$ when establishing the weak formulation. Applying $c \geqslant 1$,
the main difference would be that $1-c$ and $c$ would have different signs. In order to follow the
same strategy as in this paper one would have to find suitable weights for the space definition. In
\cite{KO1983} the authors discuss the situation when one of the weights in the weak formulation is
negative and they suggest to multiply the negative weight with another negative constant to render
it positive. Afterwards, a similar approach as in this paper could be possible.\par
In our situation there exists a second issue that may be harder to resolve. We required certain
restrictions on the growth of the function $c$, which were of the form
\begin{equation}
      \label{eq:42}
      \abs*{\partial_{z} c} \leqslant{} \kappa \sqrt{c(1-c)}
\end{equation}
for $z$ being either $x$ or $y$. The left-hand side of this inequality will always be a non-negative
real number. However, the right-hand side becomes complex-valued once $c$ exceeds 1. These growth
restrictions were important to show the coercivity of the bilinear form.\par
To conclude this section we remark that an alternative approach by means of the Helmholtz equation
already exists for the case $c>1$, see~\cite{HHBK2017}. However, this approach uses different
assumptions and yields a well-posedness theory in different spaces.
%
\subsection{Summary: What is needed to assert the existence of a solution?}
\label{sec:summary:-what-needed}
%
In this section, we summarise the necessary conditions that we had to impose on our data throughout
the paper.
\begin{enumerate}
\item The function $c$ must have weak derivatives of first order and map the domain $\Omega$ to the
      interval $[0,1]$.
\item The function $c$ must fulfil~\cref{eq:24}.
\item The function $c$ must fulfil~\cref{eq:27}.
\item The sequence $(A_{k})_{k}$ defined by~\cref{eq:31} must fulfil
      $\lim_{k\to\infty} A_{k} < 1$.
\end{enumerate}
The first three requirements can easily be verified if a concrete instance of $c$ is given. However,
the last requirement can probably only be checked in particular cases.
%
\section{Conclusion}
\label{sec:conclusion}
We have shown that a solution to the inpainting problem with the weighted Laplacian exists if the
weight is a function that maps into the interval $[0,1]$. The effort to assert the existence and
uniqueness of such a solution was significant. The well-posedness of the task can be asserted if
certain regularity conditions on the weight function $c$ are met. These requirements are similar to
what is needed to show existence and uniqueness of a solution in a discrete setting. The results in
this manuscript complete the analysis of the inpainting problem with the Laplacian. While the theory
for the discrete setup was complete for any choice of $c \geqslant 0$, the continuous theory only
covered the setup where $c > 1$. This work complements the setup where $c$ maps to $[0,1]$.
\bibliographystyle{siamplain}
\bibliography{main}
\end{document}